\documentclass[12pt]{amsart}
\usepackage[mathscr]{eucal}
\usepackage{amssymb,amsmath}
\usepackage{amsfonts}
\usepackage{a4}

\newtheorem{theorem}{Theorem}

\newtheorem{lemma}{Lemma}

\begin{document}
\title{Computing Gromov-Witten invariants of some Fano varieties.}

\author[Tomasz Maszczyk]{Tomasz Maszczyk\dag}
\address{Institute of Mathematics\\
Polish Academy of Sciences\\
Sniadeckich 8\newline 00--956 Warszawa, Poland\\
\newline Institute of Mathematics\\
University of Warsaw\\ Banacha 2\newline 02--097 Warszawa, Poland}
\email{maszczyk@mimuw.edu.pl}

\thanks{\dag The author was partially supported by KBN grants 1P03A 036 26 and 115/E-343/SPB/6.PR UE/DIE 50/2005-2008.}

\begin{abstract}
We present a recursive algorithm computing all the genus-zero
Gromov-Witten invariants from a finite number of initial ones, for
Fano varieties with generically tame semi-simple quantum (and
small quantum) $(p,p)$-type cohomology, whose first Chern class is
a strictly positive combination of effective integral basic
divisors.
\end{abstract}

\maketitle

\paragraph{\textbf{1. Introduction}} The recursive formula of Kontsevich  for numbers $N_{d}$ of rational
planar curves of degree $d$ passing through $3d-1$ points in
general position
$$N_{d}=\sum_{d_{1}+d_{2}=d}d_{1}^{2}d_{2}\left[ d_{2}\left( \begin{array}{c}
                                    3d-4 \\ 3d_{1}-2
                                    \end{array}
                                                \right) -
                           d_{1}\left( \begin{array}{c}
                                    3d-4 \\ 3d_{1}-1
                                    \end{array}
                                                \right) \right]N_{d_{1}}N_{d_{2}}$$
for $d>1$ ($N_{1}=1$ is obvious), as well as its generalizations
for genus zero Gromov-Witten invariants of all projective spaces
and Del Pezzo surfaces \cite{KoMa}, was originally derived from
the associativity relations and the homogeneity condition for
quantum cohomology \cite{KoMa} of these varieties. (For particular
Fano threefolds all Gromov-Witten invariants can be computeted
from a finite number of them by other methods \cite{Przy}.) If the
$(p,p)$-part of Hodge cohomology of a Fano variety is generated by
$H^{1,1}$ then the First Reconstruction Theorem of \cite{KoMa}
gives an inductive procedure reconstructing all restricted genus
zero Gromov-Witten classes from a finite number of codimension
zero basic classes \cite{KoMa}. In many papers quantum cohomology
of other Fano varieties was studied and in many cases completely
determined \cite{AsSad}, \cite{BaMa}, \cite{Be}, \cite{Bea},
\cite{BCF}, \cite{Bu1}, \cite{Bu2}, \cite{Cio-Fo}, \cite{GivKim},
\cite{Kim1}, \cite{Kim2}, \cite{KreTam1}, \cite{KreTam2},
\cite{SiebTi}.  But the structure of the general associativity
relations is so complicated that it is not clear how to compute
Gromov-Witten invariants in general, in particular no general
recursive structure is apparent.

However, if the Hodge $(p, p)$-part of the quantum cohomology of a
Fano variety $V$ is generically semi-simple and, moreover, admits
a tame semi-simple point lying in the subspace
$H^{1,1}(V)/(H^{1,1}(V)\cap 2\pi\sqrt{-1}H^{2}(V, \mathbb{Z}))$
(parameter space of the small quantum deformation), then all genus
zero Gromov-Witten invariants of $\bigoplus_{p}H^{p,p}(V)$ can be
reconstructed from a finite number of correlators (Reconstruction
Theorem of \cite{BaMa}).

In fact it is a consequence of the two properties: 1) the quantum
(and small quantum) $(p,p)$-type Hodge cohomology algebra is
generically tame semi-simple, 2) the system of partial
differential equations encoding associativity of quantum
multiplication on $(p,p)$-type Hodge cohomology together with the
quasi-homogeneity condition, is formally integrable. These two
properties imply that generically there exist local orthogonal
coordinates $(u_{0}, ... u_{\sigma})$ (canonical coordinates of
Dubrovin, whose basis tangent vectors form a complete system of
orthogonal idempotents), in which flatness of the metric (the
Darboux-Egoroff system of partial differential equations) together
with homogeneity, is equivalent to the Pfaff system
$$dv_{ik}=-\sum_{j\neq i, k}v_{ij}v_{kj}\ d\log\frac{u_{i}-u_{j}}{u_{k}-u_{j}}$$
for a skew-symmetric $(\sigma+1)\times (\sigma+1)$ matrix
$v=(v_{ij}=-v_{ji})_{i,j=0}^{\sigma}$. From theoretical point of
view, every diagonalizable solution $v$ reconstructs the Frobenius
manifold structure encoding the quantum cohomology up to fixing a
finite number of parameters \cite{Dub1}. The {\em isomonodromic
deformations method} based on isomonodromicity of these equations
can be found in \cite{Dub1}, \cite{Dub2}, \cite{Guzz},
\cite{Sab1}, \cite{Sab2}. However, already for $\sigma=2$ this
system reduces to the Painlev\'{e} VI equation and the inversion
formula contains inversion of complicated transcendental functions
derived from Painlev\'{e} transcendents \cite{Guzz}. In particular,
this method computes Gromov-Witten numbers $N_{d}$ of
$\mathbb{P}^{2}$ only term by term through successive expansions
inverting a complicated series \cite{Guzz}. In general, such a
series is even not known.

\vspace{3mm} In the present paper we show that the semi-simple
associativity together with the quasi-homogeneity condition reduce
to another Pfaff system
$$dy_{ab}=\sum_{c}r_{abc}(x,y)\ dx_{c}$$
with some rational functions $r_{abc}(x,y)$ with constant
coefficients, fully symmetric in indices $abc $, of a symmetric
$\sigma\times \sigma$ matrix $y=(y_{ab})_{a,b=1}^{\sigma}$
satisfying some algebraic constraint. Here independent variables
$(x_{1}, ..., x_{\sigma})$ form a part of usual affine coordinates
$(x_{0}, ..., x_{\sigma})$ on the space of the sum of
$(p,p)$-Hodge cohomology. In this way we avoid the previous hard
inversion problem from \cite{Guzz}.

In the case of $\mathbb{P}^{2}$ we solve our equation explicitly
re-obtaining the Kontsevich recursive formula for numbers $N_{d}$.
This shows that our equation does the same job as the
associativity and quasi-homogeneity equations do.

In general, using the initial data at a tame semi-simple point in
the small quantum deformation and the Newton method computing the
Taylor series of the solution, we get an algorithm computing all
Gromov-Witten numbers from a finite number of initial ones.

Our recursive algorithm generalizes the Kontsevich formula for the
projective plane and makes effective the Bayer-Manin
Reconstruction Theorem \cite{BaMa}, provided the first Chern class
of a given Fano manifold is a strictly positive combination of
effective integral basic divisors.

\vspace{3mm}
\paragraph{\textbf{2. Gromov-Witten numbers.}} Let $V$ be a  connected complex Fano manifold of dimension $n$.
Consider the complex linear space
$$H(V)=\bigoplus_{p=0}^{d}H^{p,p}(V)$$
of the $(p,p)$-Hodge cohomology algebra with the cup product $\cup
$ and with the unit $1\in H^{0,0}(V)$. We have the embedding of
additive groups
$${\rm Pic}(V)\cong H^{1,1}(V)\cap H^{2}(V, \mathbb{Z})\subset H(V)$$
(${\rm Pic}(V)\cong H^{2}(V, \mathbb{Z})$ is torsion free) and in
$H^{1, 1}(V)$ there is the cone $H_{+}^{1, 1}(V)$ of K\"{a}hler
classes. The first Chern class $c_{1}(V)$ lies in $H_{+}^{1,
1}(V)\cap H^{2}(V, \mathbb{Z})$.

According to \cite{KoMa}, given a smooth connected complex Fano
variety $V$, one chooses a basis $(h_{0},\ldots , h_{\sigma})$ of
$H(V)$ of classes Poincar\'{e} dual to classes of integral cycles
$(Z_{0},\ldots , Z_{\sigma})$ in general position, not necessarily
algebraic, such that $Z_{0}$ is the fundamental cycle,
$(Z_{1},\ldots , Z_{\rho})$ are integral effective divisors
generating the Picard group and $Z_{\sigma}$ is a point. Then
$h_{0}=1$ and $h_{1},\ldots , h_{\rho}\in H^{1, 1}(V)\cap H^{2}(V,
\mathbb{Z})$ span $H^{1, 1}(V)$. One defines the symmetric integer
valued non-degenerate symmetric matrix of the intersection form
$$g_{\alpha\beta}:=\int_{V}h_{\alpha}\cup h_{\beta},$$
and the integer valued inverse matrix $\bar{g}_{\alpha\beta}$. We
will use also fully symmetric integer valued symbols
$$g_{\alpha\beta\gamma}:=\int_{V} h_{\alpha}\cup h_{\beta}\cup h_{\gamma},$$
integers $c_{\alpha}$ such that
$$c_{1}(V)=\sum_{\alpha}c_{\alpha}h_{\alpha},$$
and the integer valued symmetric matrix
$$c_{\alpha\beta}:=\int_{V}c_{1}(V)\cup h_{\alpha}\cup h_{\beta}=\sum_{\gamma}g_{\alpha\beta\gamma}c_{\gamma}.$$
We will consider the following additional condition on the above
basis $h_{1},\ldots , h_{\rho}$ of Pic($V$).

\vspace{3mm}
\paragraph{\textbf{Condition C}} The integers $c_{\alpha}$ are
strictly positive for $\alpha=1, \ldots, \rho$.

\vspace{3mm} Such a basis exists at least for the following
classes of  Fano manifolds:

i) Generalized flag varieties $V=G/B$, where $G$ is semisimple and
$B$ is the Borel subgroup \cite{BoHi}. As $h_{1},\ldots ,
h_{\rho}$ one can take first Chern classes of homogeneous line
bundles associated with fundamental weights. Then $c_{1}(V)$
corresponds to a dominant weight $\mu_{V}$ and $c_{\alpha}=\langle
\mu_{V}, \alpha\rangle =2$, where $\alpha$ is a simple root of
$G$, which can be computed using a formula from \cite{Sno}.

ii) Fano manifolds with $\rho=1$. As $h_{1}$ one can take the
cohomology class of the ample generator of Pic($V$) and $c_{1}$ is
equal to the index of the Fano manifold $V$ .

iii) Fano threefolds with $\rho=2$. Then $V$ admits two extremal
contractions $f_{\alpha}:V\rightarrow V_{\alpha}$, corresponding
to extremal rays of lengths $\lambda_{\alpha}$, $\alpha=1, 2$,
onto projective varieties $V_{\alpha}$ with
Pic($V_{\alpha})\cong\mathbb{Z}$ and ${\rm Pic}(V)\cong
f^{*}_{1}{\rm Pic}(V_{1})\oplus f^{*}_{2}{\rm Pic}(V_{2})$. As
$h_{1}, h_{2}$ one can take pull-backs of ample generators of
Pic($V_{\alpha})$  and then $c_{1}=\lambda_{2}, c_{2}=\lambda_{1}$
\cite{MoMu}.

iv) Fano $n$-fold
$V=\mathbb{P}(\mathcal{O}_{\mathbb{P}^{1}}(1)\oplus
\mathcal{O}_{\mathbb{P}^{1}}^{n-1})$, $n\geq 2$. As a basis of
${\rm Pic}(V)$ one can take $h_{1}$ equal to the pull-back of
$\mathcal{O}_{\mathbb{P}^{1}}(1)$ and $h_{2}$ equal to the
Grothendieck tautological line bundle
$\mathcal{O}_{V/\mathbb{P}^{1}}(1)$. Then $c_{1}=1, c_{2}=n$.

v) For a toric Fano manifold $V$, corresponding to a nonsingular
complete fan with edges containing primitive vectors
$v_{1},\ldots, v_{n+\rho}$ in  $\mathbb{Z}^{n}$, one has the split
exact sequence \cite{Ful}
$$0\rightarrow \mathbb{Z}^{n}\rightarrow\mathbb{Z}^{n+\rho}\rightarrow {\rm Pic}(V)\rightarrow 0,$$
where the map $\mathbb{Z}^{n}\rightarrow\mathbb{Z}^{n+\rho}$ is
given by the matrix $V=(v_{1},\ldots, v_{n+\rho})$.

Therefore there exists an integer valued matrix $W=(w_{\alpha
i})$, $\alpha=1,\ldots,\rho$, $i=1,\ldots,n+\rho$, such that
$${\rm det}\left( \begin{array}{c}
            V \\
            W
            \end{array}\right)=\pm 1,$$
which defines integers $( c_{1}, \ldots, c_{\rho})$ as follows

$$(1, \ldots, 1)\left( \begin{array}{c}
            V \\
            W
            \end{array}\right)^{-1}=(\ast, \ldots, \ast, c_{1}, \ldots, c_{\rho}).$$
Let $D_{1},\ldots,D_{n+\rho}$ be irreducible effective invariant
divisors corresponding to the edges containing primitive vectors
$v_{1},\ldots, v_{n+\rho}$. Then divisors
$$h_{\alpha}:= \sum_{i=1}^{n+\rho}w_{\alpha i}D_{i}, $$
for $\alpha=1,\ldots,\rho$, form a basis in ${\rm Pic}(X)$ and
$$-K_{X}=\sum_{\alpha =1}^{\rho}c_{\alpha }h_{\alpha}.$$
If $W$ can be chosen in such a way that all $w_{\alpha i}\geq 0$
and $c_{\alpha }>0$ then the above condition C is satisfied.

It is easy to check that it is so for Hirzebruch  toric Del Pezzo
surfaces $\mathbb{F}_{k}$, $k=$ 1, 6, 7, where
$$v_{1}=\left( \begin{array}{c}
            0 \\
            1
            \end{array}\right),\ v_{2}=\left( \begin{array}{c}
            -1 \\
            k
            \end{array}\right),\ v_{3}=\left( \begin{array}{c}
            0 \\
            -1
            \end{array}\right),\ v_{4}=\left( \begin{array}{c}
            1 \\
            0
            \end{array}\right).$$ Namely, for
$\mathbb{F}_{1}$ (resp. $\mathbb{F}_{6}$ or $\mathbb{F}_{7}$) we
can take
$$h_{1}=D_{1}+D_{2},\ \ h_{2}=D_{1}+2D_{2}$$
(resp. $h_{1}=D_{1}$, $h_{2}=D_{2}$ or $h_{1}=D_{1}+4D_{2}$,
$h_{2}=D_{1}+5D_{2})$, and then
$$-K_{X}=h_{1}+h_{2}$$
(resp. $-K_{X}=2h_{1}+8h_{2}$  or $-K_{X}=h_{1}+h_{2}$).

\vspace{3mm}
\paragraph{\textbf{Definition.}} One
defines the number $N_{k_{1}\ldots k_{\sigma}}$ as the number of
rational curves intersecting cycles $(Z_{1},\ldots , Z_{\sigma})$
with finite multiplicities $(k_{1},\ldots, k_{\sigma})$.

One calls $N_{k_{1}\ldots k_{\sigma}}$'s the \textit{Gromov-Witten
numbers}.

\vspace{3mm}

Correctness of this definition follows from the intersection
theory on moduli stacks of stable maps \cite{BeMa}, \cite{KoMa}.
Gromov-Witten numbers are derived from the (genus zero)
\textit{Gromov-Witten invariants}. The Fano condition implies that
multiplicities $(k_{1},\ldots, k_{\rho})$ determine the homology
class of a rational curve uniquely. The collection of
Gromov-Witten invariants forms a quite complicated combinatorial
structure. Luckily, the combinatorial identities can be encoded in
properties of a generating function.

\vspace{3mm}
\paragraph{\textbf{3. Generating function.}} Let $(x_{0}, ...,
x_{\sigma})$ be the coordinate system on $H(V)$ dual to the basis
$(h_{0}, ..., h_{\sigma})$. Then the variable point of $H(V)$ has
the form $x=\sum_{\alpha}x_{\alpha}h_{\alpha}$ and one defines a
\textit{generating function} as the formal series of the form
$$F(x)=\frac{1}{6}\int_{V}x^{\cup 3} + f(x),$$
$$f(x)=\sum_{k_{1},\ldots,k_{\sigma}}N_{k_{1}\ldots k_{\sigma}}e^{k_{1}x_{1}+\ldots +k_{\rho}x_{\rho}}\frac{x_{\rho+1}^{k_{\rho+1}}}{k_{\rho+1}!}\ldots
\frac{x_{\sigma}^{k_{\sigma}}}{k_{\sigma}!},
$$
with the following condition for non-vanishing summands
$$\sum_{\alpha=1}^{\rho}k_{\alpha}c_{\alpha}+\sum_{\alpha=\rho +1}^{\sigma}k_{\alpha}(1-p_{\alpha})=3-n.$$

\vspace{3mm}
\paragraph{\textbf{4. Homogeneity.}} One defines an action of
$\mathbb{C}^{*}$ on the quotient space
$$X:=H(V)/(H^{1,
1}(V)\cap 2\pi\sqrt{-1}H^{2}(V,\mathbb{Z})) \cong
(\mathbb{C}^{*})^{\rho}\times \mathbb{C}^{\sigma+1-\rho}$$
(Dubrovin's flow) as follows
$$s\cdot \sum_{p}x_{(p)} := (x_{(1)}+\log s\ \ c_{1}(V))+\sum_{p\neq 1}s^{1-p}\ x_{(p)},$$
where $s\in \mathbb{C}^{*}$, $x_{(p)}\in H^{p,p}(V)$. The part
$f(x)$ of the generating function is homogeneous of weight $3-n$
with respect to this action \cite{KoMa}. Using the Euler field
generating the above flow
\begin{align}
E=\sum_{\alpha}E_{\alpha}\frac{\partial}{\partial x_{\alpha}}:=
\sum_{\alpha}(c_{\alpha}+(1-
p_{\alpha})x_{\alpha})\frac{\partial}{\partial x_{\alpha}}
\end{align}
one gets equivalently
$$E(f)=(3-n)f.$$
On the other hand, decomposing $x$ into homogeneous components we
get
$$E(\int_{V}x^{\cup 3}) = 3\int_{V}c_{1}(V)\cup x^{\cup 2}+(3-n)\int_{V}x^{\cup 3},$$
implying the following quasi-homogeneity property for the
generating function $F$
$$E(F) = \frac{1}{2}\int_{V}c_{1}(V)\cup x^{\cup 2}+(3-n)F,$$
which in coordinates takes the form
\begin{align}
\sum_{\delta}E_{\delta}\frac{\partial F}{\partial
x_{\delta}}=\frac{1}{2}\sum_{\alpha\beta}c_{\alpha\beta}
x_{\alpha}x_{\beta}+(3-n)F.
\end{align}
Taking second partial derivatives of this and using (1) we get the
following useful formula
\begin{align}
\sum_{\delta}E_{\delta}\frac{\partial^{3} F}{\partial
x_{\alpha}\partial x_{\beta}\partial
x_{\delta}}=c_{\alpha\beta}+(1-n+p_{\alpha}+p_{\beta})\frac{\partial^{2}
F}{\partial x_{\alpha}\partial x_{\beta}}.
\end{align}
Note that on the right hand side $c_{\alpha\beta}$ can be non-zero
only if $1-n+p_{\alpha}+p_{\beta}=0$.

\vspace{3mm}
\paragraph{\textbf{5. Flat metric.}}
The intersection form defines on $X$ a flat riemannian metric
$$g(\frac{\partial}{\partial
x_{\alpha}}, \frac{\partial}{\partial
x_{\beta}})=g_{\alpha\beta}.$$ The metric tensor is homogeneous of
weight $(2-n)$ under the Dubrovin flow.

\vspace{3mm}
\paragraph{\textbf{6. Associativity equations.}} One introduces the
following multiplication on the (trivial) tangent bundle of the
quotient space $X$
$$\frac{\partial}{\partial
x_{\alpha}}\cdot\frac{\partial}{\partial
x_{\beta}}=\sum_{\gamma,\delta}\bar{g}_{\gamma\delta}\frac{\partial^{3}
F}{\partial x_{\alpha}\partial x_{\beta}\partial
x_{\delta}}\frac{\partial}{\partial x_{\gamma}}.$$ It is
commutative by symmetry of partial derivatives and associative
with the unit $\frac{\partial}{\partial x_{0}}$ by axioms of
Gromov-Witten invariants \cite{KoMa}. The associativity condition
reads as a system of quadratic equations on third partial
derivatives of $F$. The multiplication, viewed as a tensor, is
homogeneous of weight $1$ under the Dubrovin flow, as well as the
unit of this multiplication.

The metric and the multiplication of tangent vectors satisfy
$$g(\frac{\partial}{\partial
x_{\alpha}}\cdot\frac{\partial}{\partial x_{\beta}},
\frac{\partial}{\partial x_{\gamma}}) = g(\frac{\partial}{\partial
x_{\alpha}}, \frac{\partial}{\partial x_{\beta}}\cdot
\frac{\partial}{\partial x_{\gamma}}).$$

\vspace{3mm}
\paragraph{\textbf{7. (Tame) Semi-simplicity.}} A point of $X$ is called
\textit{semi-simple} if at this point the generating function is
convergent and the above algebra of tangent vectors is
semi-simple.

According to \cite{Dub1}, locally around every semi-simple point,
there exist orthogonal coordinates $(u_{0},\ldots, u_{\sigma} )$
with a domain $U\subset X$ in which the above multiplication of
tangent vectors and the Euler vector field take the standard form
$$\frac{\partial}{\partial
u_{i}}\cdot\frac{\partial}{\partial
u_{j}}=\delta_{ij}\frac{\partial}{\partial u_{i}},$$
$$E=\sum_{i=0}^{\sigma}u_{i}\frac{\partial}{\partial
u_{i}}.$$ In particular, coordinates $(u_{0},\ldots, u_{\sigma} )$
are eigenvalues of the multiplication by the Euler field $E$. A
semi-simple point of $X$ is called \textit{tame} if at this point
$\prod_{i\neq j}(u_{i}-u_{j})\neq 0$. This condition is
independent of the choice of  coordinates $(u_{0},\ldots,
u_{\sigma} )$.

Discussion of the (tame) semi-simplicity condition, examples and
non-examples can be found in \cite{Ba}, \cite{BaMa}, \cite{Ciol},
\cite{TiXu}.

\vspace{3mm}
\paragraph{\textbf{8. Derivation of the main formula.}}
\begin{lemma}
Let the vector $(z_{i0},\ldots,z_{i\sigma})^{\intercal}$ be a
solution to the following eigen-problem with simple eigen-values
$(a_{0},\ldots, a_{\sigma} )$
\begin{align}
\sum_{\beta,\gamma}\bar{g}_{\beta\gamma}(c_{\alpha\beta}+(1-n+p_{\alpha}+p_{\beta})Y_{\alpha\beta})z_{i\gamma}=
a_{i}z_{i\alpha}.
\end{align}
for a symmetric $(\sigma+1)\times(\sigma+1)$ matrix
$Y=(Y_{\alpha\beta}=Y_{\beta\alpha})_{\alpha,\beta=0}^{\sigma}$.
Then the expression
\begin{align}
R_{\alpha\beta\gamma}(Y)=\sum_{i}\frac{z_{i\alpha}z_{i\beta}z_{i\gamma}}
{\sum_{\delta,\epsilon}\bar{g}_{\delta\epsilon}z_{i\delta}z_{i\epsilon}z_{i0}},
\end{align}
is a rational function of $Y$.
\end{lemma}

\textit{Proof:} First we rewrite the eigen-problem (12) as
$$\sum_{\beta}(a_{\alpha\beta}-a_{i}\delta_{\alpha\beta})z_{i\beta}=0.$$
The simplicity of  eigen-values  means that for every
$i=0,\ldots,\sigma$ the rank of the matrix
$(a_{\alpha\beta}-a_{i}\delta_{\alpha\beta})$ is equal to
$\sigma$. Therefore there are $\sigma$ rows of this matrix (we can
assume that they are the first $\sigma$ rows) whose exterior
product defines the one dimensional eigen-space, according to the
canonical isomorphism of vector spaces
$$\bigwedge^{\sigma}V^{*}\otimes det(V)\stackrel{\cong}{\rightarrow}V,$$
where $\dim(V)=\sigma+1$. Thus we can replace every eigen-vector
$(z_{i0},\ldots,z_{i\sigma})^{\intercal}$ in the homogeneous
function (11)  by the vector of alternated $\sigma$-minors of the
above $\sigma$ rows of the matrix
$(a_{\alpha\beta}-a_{i}\delta_{\alpha\beta})_{\alpha,\beta=0}^{\sigma}$.
It is clear that we obtain a symmetric rational function in
variables $(a_{0},\ldots, a_{\sigma})$ with coefficients in the
field of rational functions of entries $a_{\alpha\beta}$, which
are polynomials in $Y_{\alpha\beta}$'s. In this way we see that
the right hand side of (11) is a rational function of
$Y_{\alpha\beta}$'s, because we can compute this function
expressing the above symmetric rational function of eigen-values
$(a_{0},\ldots, a_{\sigma})$ in terms of their elementary
symmetric polynomials, which are polynomials in
$(a_{\alpha\beta})_{\alpha,\beta=0}^{\sigma}$, hence polynomials
in $Y_{\alpha\beta}$'s. $\Box$

Practically, to find the rational function
$R_{\alpha\beta\gamma}(Y)$ we can use here the \verb"RootSum"
function of \textit{Mathematica}. We will do so in the example of
the projective plane in paragraph 9.

 Next, let us form the
partial differential equation
\begin{align}
\frac{\partial Y_{\alpha\beta}}{\partial
x_{\gamma}}=R_{\alpha\beta\gamma}(Y)
\end{align}
and an algebraic constraint
%\begin{align}
%R_{\alpha\beta0}(Y)=g_{\alpha\beta},
%\end{align}
\begin{align}
\sum_{\gamma}E_{\gamma}R_{\alpha\beta\gamma}(Y)=c_{\alpha\beta}+(1-n+p_{\alpha}+p_{\beta})Y_{\alpha\beta}.
\end{align}

 Since the right hand side of (6) is fully symmetric in indices
 $\alpha\beta\gamma$ then locally any solution $Y$ to the
 equation (6) is of the form
\begin{align}
Y_{\alpha\beta}=\frac{\partial^{2} F}{\partial x_{\alpha}\partial
x_{\beta}}
\end{align}
for some function $F$. Since the matrix
$(a_{\alpha\beta})_{\alpha,\beta=0}^{\sigma}$ in the proof of
Lemma 1 is self-adjoint with respect to the non-degenerate
symmetric form $\bar{g}$ and has simple eigen-values its
eigen-vectors are mutually orthogonal. Therefore by (5)
\begin{align}
R_{\alpha\beta 0}(Y)=\sum_{i}\frac{z_{i\alpha}z_{i\beta}}
{\sum_{\delta,\epsilon}\bar{g}_{\delta\epsilon}z_{i\delta}z_{i\epsilon}}=g_{\alpha\beta},
\end{align}
By (5), (6), (8) and (9) $F$ automatically satisfies the
associativity equations with the unit $\frac{\partial}{\partial
x_{0}}$. Let us define
\begin{align}
f:=F-\frac{1}{6}\int_{V}x^{\cup 3}.
\end{align}

By (6), (8) and (9) we have
\begin{align}
\frac{\partial^{2}}{\partial x_{\alpha}\partial
x_{\beta}}(\frac{\partial f}{\partial x_{0}}) =0,
\end{align}
so $f$ is independent of $x_{0}$ up to adding a polynomial of
degree two, at most quadratic in $x_{0}$ and at most linear in
other variables.

By (1) we have
\begin{align}
\frac{\partial^{2}}{\partial x_{\alpha}\partial
x_{\beta}}(E(f)-(3-n)f) =
\end{align}
$$=\sum_{\gamma}E_{\gamma}\frac{\partial^{3} F}{\partial
x_{\alpha}\partial x_{\beta}\partial
x_{\gamma}}-(c_{\alpha\beta}+(1-n+p_{\alpha}+p_{\beta})\frac{\partial^{2}
F}{\partial x_{\alpha}\partial x_{\beta}}).$$ Therefore every $f$
defined in (10) by means of a solution (8) to the equation (6),
satisfying the constraint (7), is homogeneous of weight $(3-n)$ up
to adding a polynomial of degree one.

We can focus only on $f$. For this we introduce
$$y_{ab}:=Y_{ab}-\int_{V}x\cup h_{a}\cup h_{b},$$
$$r_{abc}:=R_{abc}-g_{abc},$$
where latin indices run through $(1, \ldots, \sigma)$. Then
locally
\begin{align}
y_{ab}=\frac{\partial^{2} f}{\partial x_{a}\partial x_{b}},
\end{align}
hence by (11) and (13) we have
$$\frac{\partial{y_{ab}}}{{\partial x_{0}}}=\frac{\partial^{3}f}{\partial x_{a}\partial x_{b}\partial x_{0}}=0,$$
which means that $y_{ab}$'s depend only on variables
$x=(x_{a})_{a=1}^{\sigma}$. Then the equation (6) reduces to
\begin{align}
dy_{ab}=\sum_{c}r_{abc}(x,y)dx_{c},
\end{align}
with the functions $r_{abc}(x,y)$ rational in variables
$x=(x_{a})_{a=1}^{\sigma}$, $y=(y_{ab}=y_{ba})_{a,b=1}^{\sigma}$,
and fully symmetric in indices $abc$. The algebraic constraint (7)
reduces to
\begin{align}
\sum_{c}(c_{c}+(1-p_{c})x_{c}
)r_{abc}(x,y)=n\sum_{c}g_{abc}x_{c}+(1-n+p_{a}+p_{b})y_{ab}.
\end{align}

On the other hand, in the context of semi-simple quantum
$(p,p)$-cohomology, comparing the multiplication of tangent
vectors, the unit, the metric and the multiplication by the Euler
field in coordinates $(x_{0},\ldots, x_{\sigma} )$ and
$(u_{0},\ldots, u_{\sigma} )$ we get
\begin{align}
\sum_{\gamma,\delta}\bar{g}_{\gamma\delta}\frac{\partial^{3}
F}{\partial x_{\alpha}\partial x_{\beta}\partial
x_{\delta}}\frac{\partial u_{i}}{\partial
x_{\gamma}}=\frac{\partial u_{i}}{\partial
x_{\alpha}}\frac{\partial u_{i}}{\partial x_{\beta}},
\end{align}
\begin{align}
\frac{\partial u_{i}}{\partial x_{0}} = 1,
\end{align}
\begin{align}
\sum_{\gamma,\delta}\bar{g}_{\gamma\delta}\frac{\partial
u_{i}}{\partial x_{\gamma}}\frac{\partial u_{j}}{\partial
x_{\delta}}=\delta_{ij}\sum_{\gamma,\delta}\bar{g}_{\gamma\delta}\frac{\partial
u_{i}}{\partial x_{\gamma}}\frac{\partial u_{i}}{\partial
x_{\delta}},
\end{align}
\begin{align}
\sum_{\beta,\gamma,\delta}\bar{g}_{\beta\gamma}E_{\delta}\frac{\partial^{3}
F}{\partial x_{\alpha}\partial x_{\beta}\partial
x_{\delta}}\frac{\partial u_{i}}{\partial x_{\gamma}}=
u_{i}\frac{\partial u_{i}}{\partial x_{\beta}}.
\end{align}
Equations (16)-(18) imply
\begin{align}
\frac{\partial^{3} F}{\partial x_{\alpha}\partial
x_{\beta}\partial x_{\gamma}}=\sum_{i}\frac{\frac{\partial
u_{i}}{\partial x_{\alpha}}\frac{\partial u_{i}}{\partial
x_{\beta}}\frac{\partial u_{i}}{\partial
x_{\gamma}}}{\sum_{\delta,\epsilon}\bar{g}_{\delta\epsilon}\frac{\partial
u_{i}}{\partial x_{\delta}}\frac{\partial u_{i}}{\partial
x_{\epsilon}}\frac{\partial u_{i}}{\partial x_{0}}},
\end{align}
Using (3) and (19) we get
\begin{align}
\sum_{\beta,\gamma}\bar{g}_{\beta\gamma}(c_{\alpha\beta}+(1-n+p_{\alpha}+p_{\beta})\frac{\partial^{2}
F}{\partial x_{\alpha}\partial x_{\beta}})\frac{\partial
u_{i}}{\partial x_{\gamma}}= u_{i}\frac{\partial u_{i}}{\partial
x_{\alpha}},
\end{align}
which means that for every $i$ the vector $(\frac{\partial
u_{i}}{\partial x_{0}},\ldots, \frac{\partial u_{i}}{\partial
x_{\sigma}})^{\intercal}$ is an eigenvector with the eigenvalue
$u_{i}$ of the matrix depending rationally on second partial
derivatives $\frac{\partial^{2} F}{\partial x_{\alpha}\partial
x_{\beta}}$ and constants $\bar{g}_{\alpha\beta}$, $p_{\alpha}$,
$c_{\alpha\beta}$. By Lemma 1 the right hand side of (20) is a
rational function of $\frac{\partial^{2} F}{\partial
x_{\alpha}\partial x_{\beta}}$ depending on constants
$\bar{g}_{\alpha\beta}$, $p_{\alpha}$, $c_{\alpha\beta}$. Then (3)
defines an algebraic constraint (7) on the matrix of second
derivatives of $F$.

In this way we obtain the following theorem.

\begin{theorem} The semi-simple associativity condition on a function
$$F=\frac{1}{6}\int_{V}x^{\cup 3}+f$$ for $f$ which is
\begin{itemize}
\item  independent of $x_{0}$ up to adding a polynomial of degree two, at most quadratic in $x_{0}$ and at most linear in
other variables,
\item  homogeneous of weight $(3-n)$ up to adding a polynomial of degree
one,
\end{itemize}
is locally equivalent to the problem (14)-(15).
\end{theorem}

\vspace{3mm}
\paragraph{\textbf{9. Example $V=\mathbb{P}^{1}$.}} The basis $(h_{0},h_{1})$ in
$(p,p)$-Hodge cohomology is defined  uniquely by the above
convention. We have $n=1$, $p_{0}=0, p_{1}=1$, $g_{00}=g_{11}=0,\
g_{01}=1$, $g_{000}=g_{011}=g_{111}=0,\ g_{001}=1$, $c_{0}=0,
c_{1}=2$, and $c_{00}=2,\ c_{01}=c_{11}=0$.

Then we obtain
$$r_{111}(x,y)=y_{11},$$
so the system (14) has the form
\begin{align}
dy_{11}=y_{11}dx_{1},
\end{align}
and the constraint (15) is satisfied automatically. For
$$f(x_{1})=\sum_{k_{1}}N_{k_{1}}e^{k_{1}x_{1}}$$
the condition on non-vanishing terms is
$$2k_{1}=2,$$
which determines $f$ uniquely up to a constant factor, as well as
the system (14) does.

\vspace{3mm}
\paragraph{\textbf{10. Example $V=\mathbb{P}^{2}$.}} The basis $(h_{0}, h_{1}, h_{2})$ in
$(p,p)$-Hodge cohomology as above is defined uniquely by the
condition that the cycle $Z_{1}$ is a hyperplane. Now $n=2$,
$p_{\alpha}=\alpha$, $g_{\alpha\beta}=1$ for $\alpha+\beta=2$,
$g_{\alpha\beta\gamma}=1$ for $\alpha+\beta+\gamma=2$,
$c_{\alpha}=3$ for $\alpha=1$, $c_{\alpha\beta}=3$ for
$\alpha+\beta=1$,  and all symbols are zero otherwise.

 Then we obtain (using \textit{Mathematica})
\begin{align*}
r_{111}(x,y)& =\frac{9y_{11}+x_{2}(y_{11}^{2}+6y_{12})+3x_{2}^{2}y_{22}}{27+3x_{2}y_{11}-2x_{2}^{2}y_{12}},\\
r_{112}(x,y)& =\frac{18y_{12} +x_{2}(2y_{11}y_{12}+9y_{22})}{27+3x_{2}y_{11}-2x_{2}^{2}y_{12}},\\
r_{122}(x,y)& =\frac{27y_{22}+4x_{2}y_{12}^{2}}{27+3x_{2}y_{11}-2x_{2}^{2}y_{12}},\\
r_{222}(x,y)&
=\frac{12y_{12}^{2}-9y_{11}y_{22}+6x_{2}y_{12}y_{22}}{27+3x_{2}y_{11}-2x_{2}^{2}y_{12}},
\end{align*}
so the Pfaff system (14) has the form
\begin{align}
dy_{11}& =\frac{9y_{11}+x_{2}(y_{11}^{2}+6y_{12})+3x_{2}^{2}y_{22}}{27+3x_{2}y_{11}-2x_{2}^{2}y_{12}}dx_{1}+\frac{18y_{12} +x_{2}(2y_{11}y_{12}+9y_{22})}{27+3x_{2}y_{11}-2x_{2}^{2}y_{12}}dx_{2},\\
dy_{12}& =\frac{18y_{12}
+x_{2}(2y_{11}y_{12}+9y_{22})}{27+3x_{2}y_{11}-2x_{2}^{2}y_{12}}dx_{1}
        +\frac{27y_{22}+4x_{2}y_{12}^{2}}{27+3x_{2}y_{11}-2x_{2}^{2}y_{12}}dx_{2},\\
dy_{22}&
=\frac{27y_{22}+4x_{2}y_{12}^{2}}{27+3x_{2}y_{11}-2x_{2}^{2}y_{12}}dx_{1}+\frac{12y_{12}^{2}-9y_{11}y_{22}+6x_{2}y_{12}y_{22}}{27+3x_{2}y_{11}-2x_{2}^{2}y_{12}}dx_{2}.
\end{align}
One can check (using \textit{Mathematica}) that the algebraic
constraint (15) is satisfied identically.

The system (23)-(25) is equivalent to the system
\begin{align}
&
d(x_{2}y_{11})=\frac{9x_{2}y_{11}+x_{2}^{2}(y_{11}+6y_{12})+3x_{2}^{3}y_{22}}{27+3x_{2}y_{11}-2x_{2}^{2}y_{12}}
d(x_{1}+3\log x_{2})
,\\
& d(3x_{2}y_{11}-x_{2}^{2}y_{12})
=x_{2}y_{11}d(x_{1}+3\log x_{2}),\\
& d(18x_{2}y_{11}-9x_{2}^{2}y_{12}+x_{2}^{3}y_{22})
=2(3x_{2}y_{11}-x_{2}^{2}y_{12}) d(x_{1}+3\log x_{2}).
\end{align}

If we introduce the expression $\phi$ such that
\begin{align}
2\phi:=18x_{2}y_{11}-9x_{2}^{2}y_{12}+x_{2}^{3}y_{22},
\end{align}
then the system (26)-(28) means that
\begin{align}
\phi & =\phi(x_{1}+3\log x_{2}),\\
\phi'''&
=\frac{6\phi-33\phi'+54\phi''+\phi''^{2}}{27+2\phi'-3\phi''},
\end{align}
which can be written equivalently as
\begin{align}
27\phi'''-54\phi''+33\phi'-6\phi=3\phi''\phi'''-2\phi'\phi'''+\phi''^{2}.
\end{align}
 For
$$f(x_{1},x_{2})=\sum_{k_{1},k_{2}}N_{k_{1}k_{2}}e^{k_{1}x_{1}}\frac{x_{2}^{k_{2}}}{k_{2}!}$$
the condition on non-vanishing terms is
$$3k_{1}-k_{2}=1.$$
 Since $Z_{1}$ is a hyperplane
$k_{1}$ equals to the degree $d$ of a rational curve, hence
$k_{2}=3d-1$, $N_{k_{1}k_{2}}=N_{d}$ and
$$f=\sum_{d}N_{d}e^{dx_{1}}\frac{x_{2}^{3d-1}}{(3d-1)!}.$$
Now we compute $\phi$ substituting $f$ into (13) and next the
obtained result  to (29). We get
\begin{align}
\phi=\sum_{d}N_{d}e^{dx_{1}}\frac{x_{2}^{3d}}{(3d-1)!}=\sum_{d}\frac{N_{d}}{(3d-1)!}e^{d(x_{1}+3\log
x_{2}) }.
\end{align}
Inserting this to (32) we get finally
$$N_{d}=\sum_{d_{1}+d_{2}=d}d_{1}^{2}d_{2}\left[ d_{2}\left( \begin{array}{c}
                                    3d-4 \\ 3d_{1}-2
                                    \end{array}
                                                \right) -
                           d_{1}\left( \begin{array}{c}
                                    3d-4 \\ 3d_{1}-1
                                    \end{array}
                                                \right) \right]N_{d_{1}}N_{d_{2}}.$$
This shows that our system provides the same answer as the
original use of associativity and homogeneity conditions.

\vspace{3mm}
\paragraph{\textbf{11. Generalization.}} In
general, we can use the Newton method to compute the Taylor series
of the solution. To apply this we have to know that the
denominator of the rational function on the right hand side of our
equation (15) do not vanish. It is so at a tame semi-simple point
because all denominators on the right hand side of (20) don't
vanish by (17), (18) and non-degeneracy of the matrix $\bar{g}$.
Computing derivatives of the solution $y$ which are derivatives of
the generating function $f$, at a tame semi-simple point
$(x_{1},\ldots, x_{\rho}, 0,\ldots, 0)$ in the parameter space of
the  small quantum deformation
\begin{align}
\frac{\partial^{m_{1}+\ldots +m_{\sigma}}f}{\partial
x_{1}^{m_{1}}\ldots\partial
x_{\sigma}^{m_{\sigma}}}\mid_{x_{\rho+1}=\ldots=x_{\sigma}=0}=\sum_{k_{1},\ldots,
k_{\rho} }k_{1}^{m_{1}}\ldots
k_{\rho}^{m_{\rho}}e^{k_{1}x_{1}+\ldots+
k_{\rho}x_{\rho}}N_{k_{1}\ldots k_{\rho}m_{\rho+1}\ldots
m_{\sigma}}
\end{align}
with the following condition on non-vanishing summands
\begin{align}\sum_{a=1}^{\rho}k_{a}c_{a}=\sum_{a=\rho +1}^{\sigma}m_{a}(p_{a}-1)+3-n,
\end{align}
we see, provided the condition C is satisfied, that on the right
hand side of (34) we have a finite sum, because there is only a
finite number of non-negative vectors $(k_{1},\ldots, k_{\rho})$
which are solutions to (35).

By the Lagrange interpolation we find polynomials
$$P_{l_{1}\ldots l_{\rho}}(t_{1},\ldots,t_{\rho})\in
\mathbb{Q}[t_{1},\ldots, t_{\rho}]$$ such that for all (a finite
number of) non-negative solutions $(k_{1},\ldots, k_{\rho})$ to
(35)
\begin{align}P_{l_{1}\ldots
l_{\rho}}(k_{1},\ldots,
k_{\rho})=\delta_{k_{1}l_{1}}\ldots\delta_{k_{\rho}l_{\rho}}.\end{align}
If we expand these polynomials as follows
\begin{align}P_{l_{1}\ldots
l_{\rho}}(t_{1},\ldots, t_{\rho})=\sum_{m_{1},\ldots
 m_{\rho}}a_{l_{1}\ldots l_{\rho}, m_{1}\ldots
m_{\rho}}t_{1}^{m_{1}}\ldots t_{\rho}^{m_{\rho}}\end{align} and
apply this expansion to monomials $k_{1}^{m_{1}}\ldots
k_{\rho}^{m_{\rho}}$ on the right hand side of (34) substituted
instead of $t_{1}^{m_{1}}\ldots t_{\rho}^{m_{\rho}}$, we get by
(36) the following expression for the Gromov-Witten numbers
\begin{align}
 & N_{k_{1}\ldots k_{\sigma}} =\end{align}
 \begin{align*}= e^{-(k_{1}x_{1}+\ldots+k_{\rho}x_{\rho})} &
 \sum_{m_{1},\ldots,
 m_{\rho}}a_{k_{1}\ldots k_{\rho}, m_{1}\ldots
m_{\rho}}\frac{\partial^{m_{1}+\ldots +m_{\rho}+k_{\rho+1}+\ldots
+k_{\sigma}}f}{\partial x_{1}^{m_{1}}\ldots\partial
x_{\rho}^{m_{\rho}}\partial x_{\rho+1}^{k_{\rho+1}}\ldots\partial
x_{\sigma}^{k_{\sigma}}}\mid_{x_{\rho+1}=\ldots=x_{\sigma}=0}.
\end{align*}
Therefore the process of consecutive differentiation of our
equation and substitution of already computed derivatives is a
recursive algorithm  computing all Gromov-Witten numbers from a
finite number of initial ones. By (34) and (35) the set of initial
Gromov-Witten numbers consists of
\begin{align}
& N_{k_{1}\ldots k_{\rho}0\ldots 0} & \ \ \ \ & {\rm where}\ \
\sum_{c=1}^{\rho}k_{c}c_{c}  =  3-n, \\
& N_{k_{1}\ldots k_{\rho}0\ldots 1_{a}\ldots 0} & \ \ \ \ & {\rm
where}\ \ \sum_{c=1}^{\rho}k_{c}c_{c}=
p_{a}+2-n, \\
& N_{k_{1}\ldots k_{\rho}0\ldots 1_{a}\ldots 1_{b}\ldots 0} & \ \
\ \ & {\rm where}\ \
\sum_{c=1}^{\rho}k_{c}c_{c}  =  p_{a}+p_{b}+1-n, \\
& N_{k_{1}\ldots k_{\rho}0\ldots 2_{a}\ldots 0} & \ \ \ \ & {\rm
where}\ \ \sum_{c=1}^{\rho}k_{c}c_{c}  = 2p_{a}+1-n.\end{align}

In view of paragraph 10, this can be regarded as a generalization
of the Kontsevich recursive formula for the projective plane.

\vspace{3mm}
\paragraph{\textbf{Remark.}} For a Fano $n$-fold of index $r$,
among the initial Gromov-Witten numbers the numbers
$N_{k_{1}\ldots k_{\rho}0\ldots 0}$ (resp. $N_{k_{1}\ldots
k_{\rho}0\ldots 1_{a}\ldots 0}$, $N_{k_{1}\ldots k_{\rho}0\ldots
1_{a}\ldots 1_{b}\ldots 0}$, $N_{k_{1}\ldots k_{\rho}0\ldots
2_{a}\ldots 0}$) can appear only if $n+r\leq 3$ (resp. $n+r-2\leq
p_{a}\leq n$, $n+r-1\leq p_{a}+p_{b}\leq 2n-1$, $(n+r-1)/2 \leq
p_{a}\leq n$).

In particular, for $\mathbb{P}^{n}$, $n\geq 2$, the only initial
Gromov-Witten number with respect to the basic classes $H,
H^{2},\ldots, H^{n}$ is the number $N_{10\ldots 02}=1$ of lines
passing through two points in general position.

\end{document}